\documentclass[letterpaper, 10 pt, conference]{ieeeconf}  

\IEEEoverridecommandlockouts 

\usepackage{afterpage}
\usepackage{epsfig} 
\usepackage{times} 
\usepackage{amsmath} 
\usepackage{amssymb}  
\usepackage{amsfonts}
\usepackage{mathptmx} 
\usepackage{tabularx}
\usepackage{epstopdf}
\usepackage{float}
\usepackage{epic,color}
\usepackage{mathrsfs}
\usepackage{dsfont,MnSymbol}
\usepackage{algorithm}
\usepackage{multirow} 
\usepackage[T1]{fontenc}

\newcounter{rmnum}
\newenvironment{romannum}{\begin{list}{{\upshape (\roman{rmnum})}}{\usecounter{rmnum}
\setlength{\leftmargin}{6pt}
\setlength{\rightmargin}{0pt}
\setlength{\itemsep}{1pt}
\setlength{\itemindent}{0pt}
}}{\end{list}}

\newcounter{anum}

\def\IEEEQEDclosed{\mbox{\rule[0pt]{1.3ex}{1.3ex}}}

\def\qed{\ifmmode\IEEEQEDclosed\else{\unskip\nobreak\hfil
		\penalty50\hskip1em\null\nobreak\hfil\IEEEQEDclosed
		\parfillskip=0pt\finalhyphendemerits=0\endgraf}\fi}

\def\qed{\hspace*{\fill}~\IEEEQED\par\endtrivlist\unskip}

\def\Re{\mathbb{R}}

\def\Lemma#1{Lemma~\ref{#1}}

\def\Sec#1{Sec.~\ref{#1}}

\def\notes#1{\marginpar{\tiny #1}\typeout{Notes!
Notes!
Notes!
}}
\renewcommand{\notes}[1]{\typeout{notes!}}

\def\Re{\field{R}}

\def\Sec#1{Sec.~\ref{#1}}

\def\clP{{\cal P}}
\def\clZ{{\cal Z}}

\def\Sec#1{Sec~\ref{#1}}

\def\E{{\sf E}}

\def\clC{{\cal C}}

\def\Sec#1{Sec.~\ref{#1}}

\def\IEEEQEDclosed{\mbox{\rule[0pt]{1.3ex}{1.3ex}}}
\def\qed{\nobreak\hfill\IEEEQEDclosed}

\def\clZ{{\cal Z}}

\newtheorem{theorem}{Theorem}

\newtheorem{definition}{Definition}
\newtheorem{lemma}{Lemma}
\newtheorem{remark}{Remark}
\newtheorem{proposition}{Proposition}

\def\beq{\begin{eqnarray}} 
\def\bc{\begin{center}} 
\def\be{\begin{enumerate}}
\def\bi{\begin{itemize}} 
\def\bs{\begin{small}}
\def\bS{\begin{slide}}
\def\ec{\end{center}} 
\def\ee{\end{enumerate}}
\def\ei{\end{itemize}}
\def\es{\end{small}}
\def\eS{\end{slide}}
\def\eeq{\end{eqnarray}}

\newcommand{\newP}[1]{\medskip\noindent{\bf #1:}}

\newcommand{\ud}{\,\mathrm{d}}

\def\Re{\mathbb{R}}
\def\E{{\sf E}}

\def\Sec#1{Sec.~\ref{#1}}
\def\Thm#1{Thm.~\ref{#1}}
\def\Prop#1{Prop.~\ref{#1}}

\def\clP{{\cal P}}
\def\clZ{{\cal Z}}

\renewcommand{\Re}{\mathbb{R}}

\newcommand{\tf}{{t_f}}

\newcommand{\sfJ}{{\sf J}}

\def\sJ{{\sf J}}
\def\sP{{\sf P}}

\def\ones{{\sf 1}}

\def\clC{{\cal C}}
\def\clF{{\cal F}}

\def\clP{{\cal P}}
\def\clU{{\cal U}}

\def\clZ{{\cal Z}}
\def\E{{\sf E}}

\def\bS{\mathbb{S}}
\def\dv{\operatorname{diag}}
\def\sp{\operatorname{span}}

\def\barpi{{\bar{\mu}}}
\def\tf{{\cal T}}

\newlength{\noteWidth}
\long\def\notes#1{\ifinner
	{\tiny #1}
	\else
	\marginpar{\parbox[t]{\noteWidth}{\raggedright\tiny #1}}
	\fi}

\title{\LARGE \bf A Dual Characterization of the Stability of the Wonham Filter}

\author{Jin Won Kim and Prashant G. Mehta
	\thanks{Financial support from the 
		NSF grant 1761622 and the 
		ARO grant W911NF1810334 is gratefully acknowledged. 
	}
	\thanks{J.~W. Kim and P.~G.~Mehta are with the Coordinated
		Science Laboratory and the Department of Mechanical Science and
		Engineering at the University of Illinois at Urbana-Champaign
		(UIUC); Corresponding email: mehtapg@illinois.edu.}
}

\begin{document}

\maketitle
\thispagestyle{empty}
\pagestyle{empty}

\begin{abstract}

This paper revisits the classical question of the stability of the
nonlinear Wonham filter.  The novel contributions of this paper are two-fold:
(i) definition of the stabilizability for the (control-theoretic) dual to
the nonlinear filter; and (ii) the use of this definition to obtain conclusions on the
stability of the Wonham filter.  Specifically, it is shown that the
stabilizability of the dual system is necessary for filter stability
and conversely stabilizability implies that the filter asymptotically
detects the correct ergodic class. The formulation and the proofs are
based upon a recently discovered duality result whereby the nonlinear
filtering problem is cast as a stochastic optimal control problem for
a backward stochastic differential equation (BSDE).  The
control-theoretic proof
techniques and results
may be viewed as a generalization of the classical work on the stability
of the Kalman filter.  

	
\end{abstract}

\section{Introduction}
\label{sec:intro}

Viewed from a certain lens, the story of stochastic filter stability
is a story of two parts: (i) stability of the Kalman filter where
control-theoretic definitions and methods are paramount; and (ii)
stability of the nonlinear filter where there is no hint of such
methods (with one notable exception~\cite{van2006filtering}). 
Arguably, the control techniques are useful for 
the analysis of Kalman filter, because of the classical dual relationship
between observability and controllability of deterministic linear
systems.  This dual relationship extends to the stochastic linear Gaussian settings: 
In Kalman's celebrated paper with Bucy, it is shown that the Kalman
filter is dual to a deterministic linear quadratic (LQ) optimal control problem.
The relationship is useful in two ways: (i) Asymptotic
stability of the filter is related to the asymptotic
stability of the dual optimal control system; and (ii) Necessary and
sufficient conditions for the same are stabilizability for
the optimal control problem, and (because of the dual relationship)
detectability for the filter.  
Notably, duality explains why, with the
time arrow reversed, the covariance update equation of the Kalman
filter is the same as the dynamic Riccati equation (DRE) of optimal
control.  
In practical terms, asymptotic stability of the Kalman filter is deduced by establishing an asymptotic limit for the value
function of the dual LQ problem~\cite[Ch.~9]{xiong2008introduction}.   


Our goal in this paper is to extend these classical control theoretic
techniques for the stability analysis of the nonlinear Wonham filter.
Specifically, we are interested in obtaining necessary and sufficient
conditions for filter stability.  Our focus is on the so called
non-ergodic signal case and we are interested in a minimal condition
that a model should satisfy for the filter stability to hold.  In
the filter stability literature, this condition is referred to as 
detectability~\cite{van2009observability}.


This problem has a rich and storied history; 
cf.,~\cite{van2010nonlinear,chigansky2009intrinsic} and references
therein.  For the non-ergodic signal case, a notable early
contribution is~\cite{clark1999relative} where formulae for the
relative entropy are derived and it is shown that the relative entropy
is a Lyapunov function for the filter.  
Our paper is closely inspired by~\cite{baxendale2004asymptotic} who
are the first to formulate certain ``identifying conditions'' that
are shown to be sufficient for the stability of Wonham filter.  These
conditions are formulated in terms of the model parameters
(transition matrix and the observation function).  The definition of observability and detectability appears
in~\cite{van2009observability,van2009uniform}.  For the Wonham filter, the subspace
of observable functions is completely characterized in terms of model
parameters.  Detectability
is shown to be both necessary and sufficient for filter stability.
Extensions to these definitions have recently appeared
in~\cite{mcdonald2018stability,mcdonald2019cdc}.

The paper has a single contribution given as~\Thm{thm:main-result}:
We define the stabilizability property of the dual system and relate
this property to the asymptotic stability of the
filter. Stabilizability is shown to be the dual to the detectability
definition of~\cite{van2009observability}.  
The overall development -- introduction of the dual
system, stabilizability definition, and its use the in the filter stability analysis -- has close parallels to the Kalman
filter stability theory.  This connection is explained using several remarks
in the paper.  While the narrow focus of this paper is on the 
non-ergodic signal case (where stabilizability is non-trivial), a companion paper presents filter stability results
for the ergodic signal case~\cite{kim2021ergodic}. 
The analysis in both these papers is based upon a recently discovered
duality result whereby the nonlinear filtering problem is cast as a
stochastic optimal control problem for a backward stochastic
differential equation (BSDE)~\cite{kim2019duality}.   


Both the optimal control formulation and its use in obtaining the
filter stability proofs are new.  While~\cite{van2006filtering}
also employed a dual optimal control problem, it is completely
different from the dual formulation used here.  As explained in
our earlier papers~\cite{kim2019duality,kim2020smoothing,kim2019observability},
our formulation is a generalization of the Kalman-Bucy
duality while~\cite{van2006filtering} is a generalization
of the minimum energy or maximum likelihood duality (see
also~\cite{todorov2008general}).       
Notably, the
classical proofs of the stability of the Kalman filter are based on
the original Kalman-Bucy duality (see e.g.~~\cite{xiong2008introduction}).  
Our proofs can thus be viewed as
a generalization of the linear stability theory.  



The outline of the remainder of this paper is as follows: The problem
formulation appears in \Sec{sec:prelim}. 
The background on duality and stabilizability
is in~\Sec{sec:duality}.   The main result and its proof
are in \Sec{sec:main} and \Sec{sec:main_proof}, respectively.

\section{Problem formulation}\label{sec:prelim}

\newP{Notation} The state-space $\bS:=\{1,2,\hdots,d\}$ is finite.  
The set of probability vectors on $\bS$ is denoted by $\clP(\bS)$:
$\mu\in \clP(\bS)$ if $\mu(x)\geq 0$ and $\sum_{x\in\bS} \mu(x) =
1$. The space of deterministic functions on $\bS$ is identified with
$\Re^d$: Any function $f:\mathbb{S}\to \Re$ is determined by its
value $f(x)$ at $x\in \bS$.  For a measure $\mu \in {\cal
  P}(\mathbb{S})$ and a function $f\in \Re^d$,   $\mu(f):=\sum_x
\mu(x) f(x)$. For two vectors $f,h\in\Re^d$, $f\, h$ denotes the
element-wise (Hadamard) product: $(f\,h)(x) := f(x) h(x)$ and
similarly $f^2  = f\, f$.  
The vector of all ones is denoted as $\ones$ and $f\mid_{\ones^\perp}:=f
- \frac{1}{d} (1^\top f) 1$.  For a subset
${D} \subset \mathbb{S}$, $\ones_{{D}}$ denotes the
indicator function with support on ${D}$.

\subsection{Filtering model}\label{ssec:model}

Consider a pair of continuous-time stochastic processes $(X,Z)$
defined on a probability space $(\Omega,\clF,{\sf P})$:

(1) The state $X= \{X_t\in\bS:t\geq 0\}$ is a Markov process with
  initial condition $X_0\sim\mu\in{\clP}(\bS)$ (prior) and the
  generator (row stochastic rate matrix) 
  $A\in\Re^{d\times d}$.  The finite state-space 
  is partitioned into $m$ ergodic classes $\bS = \bigcup_{k=1}^m
  \bS_k$ such that $\sP([X_t\in \bS_l] \mid [X_0\in \bS_k]) = 0$ for
  all $t\ge 0$ and $l \neq k$. 
For a function $f\in \Re^d$, the carr\'e du
  champ operator $\Gamma:\Re^d\to\Re^d$ is defined according to 
\[\Gamma (f)(x) := \sum_{j \in \mathbb{S}} A(x,j)
  (f(x) - f(j))^2 \quad \text{for} \; x\in\bS 
\] 

(2) The observation process $Z = \{Z_t\in\Re:t\geq 0\}$ is defined
  according to the following model:
\[
Z_t := \int_0^t h(X_s) \ud s + W_t
\]
where $h : \bS \to \Re$ is the observation function and 
$W=\{W_t\in\Re:t\geq 0\}$ is a Wiener process (w.p.) that is assumed
to be independent of $X$. The
  scalar-valued observation model is considered for notational
  ease.  The covariance of $W$ is denoted as $R$ which is assumed to
be positive.  The filtration
generated by $Z$ is denoted as $\clZ := \{\clZ_t : 0\le t\le T\}$
where $\clZ_t = \sigma(\{Z_s: 0\le s\le t\})$.


\newP{Function spaces} To stress the prior $\mu$, the
probability measure is denoted $\sP^\mu$ and the 
expectation is ${\sf E}^{\mu}(\cdot)$. The
space of square-integrable $\clZ_T$-measurable random
functions on $\bS$ is denoted $L^2_{\clZ_T}(\bS)$:  $F\in
L^2_{\clZ_T} (\bS)$ if $F$ is $\clZ_T$-measurable and ${\sf
  E}^{\mu}(|F(X_T)|^2) < \infty$.  The space of $\clZ$-adapted
square-integrable $S$-valued stochastic processes is denoted 
$L^2_{\clZ}([0,T];S)$.  Typical examples of $S$ are (i) $\Re$ for
real-valued, and (ii) $\Re^d$ for vector-valued stochastic processes. 

\medskip

The filtering problem is to compute the conditional 
distribution (posterior), denoted
$\pi_t^{\mu}\in\clP(\bS)$, of the state $X_t$ given $\clZ_t$.  For $f\in  \Re^d$, $
\pi_t^{\mu}(f) := \E^{\mu}(f(X_t)|\clZ_t)
$ is the object of interest.

\def\lsq#1#2{L_{#1}^2([0,T]\,;{#2})}

\medskip

\subsection{Definition of filter stability}

In the finite state-space settings of this paper, the Wonham filter is
a stochastic differential equation (SDE):
\begin{align}
\ud \pi_t &= A^\top \pi_t \ud t + \big(\dv(h)-\pi_t h^\top
            \big)\,\pi_t\, R^{-1}\,\big(\ud Z_t- \pi_t(h) \ud t\big)\label{eq:Wonham}
\end{align}
With an initialization $\pi_0=\nu  \in \clP(\bS)$, the solution of the
Wonham filter is denoted as $\pi^\nu := \{\pi_t^\nu \in\clP(\bS):t\geq 0\}$.  The
(true) posterior $\pi^{\mu}:= \{\pi_t^\mu \in\clP(\bS):t\geq 0\}$ results from the choice of
the initial condition $\pi_0 = \mu$.

\medskip

\begin{definition}\label{def:filter-stability-l2}
The Wonham filter is {\em stable}, in the sense of weak convergence in $L^2$, if for each $f \in \Re^d$
\begin{equation}\label{eq:filter-stability}
\E^\mu \big(| \pi_t^\mu(f) - \pi_t^\nu(f) |^2\big) \longrightarrow 0 \quad \text{as }t\to\infty
\end{equation}
whenever $\mu \ll \nu$.
\end{definition}

\medskip


\section{Duality and stabilizability}\label{sec:duality}

\subsection{Background on duality for nonlinear filtering}

In our recent work~\cite{kim2019duality}, a dual optimal control
problem is introduced.  It is based upon the following
backward stochastic differential equation (BSDE):
\begin{align}\
- \ud {Y}_t(x)    & = \left((A {Y}_t)(x) + h(x)U_t + h (x) 
          {V}_t(x)\right)\ud t - {V}_t(x)\ud Z_t, 
\nonumber \\
Y_T (x) &= f(x) \;\; \forall \; x \in \bS \label{eq:dual-bsde}
\end{align}
where 
the control input $U=\{U_t\in \Re\,:\, 0\le t \le T\}$ is
in $L^2_{\clZ}([0,T];\Re)=: {\cal U}$.  Such controls are referred to as
{\em admissible}. 
The solution $(Y,V):=\{(Y_t,V_t) \in \Re^d\times \Re^d\,:\, 0\le t \le
T\}$ is in $L^2_{\clZ}([0,T];\Re^d\times \Re^d)$.  That is, the
solution is forward-adapted to the filtration $\clZ$.

Consider next the following estimator (for the random variable $f(X_T)$):
\begin{equation}\label{eq:opt_est}
S_T = \pi_0(Y_0) - \int_0^T U_t \ud Z_t
\end{equation}
The estimator is parameterized with an admissible control
$U \in {\cal U}$ and a given $\pi_0\in {\cal
  P}(\mathbb{S})$.  $Y_0$ is obtained by
solving the BSDE~\eqref{eq:dual-bsde} with the control input $U$ and $Y_T=f$.  

By an application of the It\^{o}-Wentzell theorem (see~\cite[Appendix~A]{kim2019duality})
$$
f(X_T) - S_T = Y_0(X_0) - \pi_0(Y_0) + \int_0^T (U_t + V_t(X_t)) \ud W_t + Y_t^\top \ud N_t
$$
where 
$N$ is a ${\sf P}^\mu$-martingale. 
Upon squaring and taking expectation, we obtain the duality relationship
\begin{align}\label{eq:duality_reln}
\E^{\mu} ( |f(X_T) &- S_T|^2 )= \sJ_T^\mu (U;f) +
                             |\pi_0(Y_0) -\mu(Y_0)|^2
\end{align}
where
\begin{align*}
\sJ_T^\mu (U;f)
:= \E^\mu \Big( |Y_0(X_0)-\mu(Y_0)|^2 + \int_0^T \ell
  (Y_t,V_t,U_t\,;X_t) \ud t \Big)
\end{align*}
where $\ell(y,v,u\,;x) = \Gamma(y)(x) + |u+v(x)|_R^2$;
cf.,~\cite[Eq.~(6)]{kim2019duality}. 
In~\cite{kim2019duality},~\eqref{eq:duality_reln} is referred to
as the {\em duality principle}.  


\medskip

The {\em dual optimal control problem} is to choose a 
control $U\in{\cal U}$ such that $\sJ_T^\mu (U;f)$ is minimized subject to the
BSDE constraint~\eqref{eq:dual-bsde}.  The minimum value is denoted as
$\sJ_T^\mu(f)$ or more simply as $\sJ_T^\mu$ if there is no
chance of ambiguity.  

The existence and uniqueness of the optimal control follows from the
standard results in the BSDE constrained optimal control
theory~\cite{peng1993backward}.  The solution, including the formula
for optimal control, is described in~\cite[Theorem~1]{kim2019duality}.
Some technical background is also included in the Appendix~\ref{apdx:dual-optimal-control} of
this paper.  
Let $U^{\mu} :=\{U_t^{\mu}:\, 0\le t \le T\}$ be the
optimal control and $(Y^\mu,
V^\mu):=\{({Y}_t^\mu,{V}_t^\mu):\,
0\le t \le T\}$ is the associated (optimal) trajectory.  Then:

(1)~\cite[Theorem~2]{kim2019duality}: The optimal control gives the conditional mean
\begin{equation}\label{eq:optimal-estimate-mu}
\pi_T^\mu(f) = \mu(Y_0^\mu) - \int_0^T U_t^{\mu} \ud Z_t,\quad {\sf P}^\mu-\text{a.s.}
\end{equation} 

(2)~\cite[Theorem~5]{kim2019duality}: The optimal value function
  gives the variance
\[
\sJ_T^\mu(f)  = \E^{\mu}( |f(X_T) - \pi_T^\mu(f) |^2 )
\] 
Using this formula, the optimal value is uniformly bounded
by $\frac{1}{4}|\text{osc}(f)|^2$ where 
$\text{osc}(f):=\max_{i,j\in\bS}|f(i)-f(j)|$. 

\medskip


\subsection{Stabilizability of the BSDE}
\label{sec:obsvbl}

\begin{definition}\label{def:controllable-subspace}
For the BSDE~\eqref{eq:dual-bsde}, the 
\emph{controllable subspace} is defined as:
\begin{equation*}
\clC_T :=\{y_0\in\Re^d \;\mid \; \exists \; c \in \Re,\, U \in {\cal U} \;\;
            \text{s.t.} \;\; Y_0 = y_0
   \; \; \text{and} \;\; Y_T=c\ones \}
\end{equation*}
(Note $c,y_0$ are deterministic and $U$ is an admissible stochastic process.) If $\clC_T =\Re^d$ then the BSDE is said to be {\em controllable}.  
\end{definition}

\medskip

Because $1\in\clC_T$, it is a non-trivial subspace of $\Re^d$
(even with $h=0$). 
An explicit characterization of the controllable subspace is given in
the following: 

\medskip

\begin{proposition}{\cite[Theorem~2]{kim2019observability}} 
\label{thm:rank}
For any positive value of terminal time $T$, $\clC_T$ is the smallest such subspace
$\clC\subseteq \Re^d$ that satisfies the following two properties:
\begin{romannum}
	\item The constant function $\ones\in \clC$; and
	\item If $g\in\clC$ then $A g \in \clC$ and $g h
          \in\clC$. 
\end{romannum}
Explicitly,
	\begin{align*}
	\clC  := \sp\big\{\ones, &\,  h, \,  Ah, \,  A^2h, \,  A^3h, \, \ldots, \label{eq:obs_gram_nl}\\
	&h^2, \,  A(h^2), \,  h (Ah), \,  A^2(h^2),\ldots, \nonumber\\
	&h^3, \,  (Ah) (h^2), \,  h A(h^2), \, \ldots \big\} \nonumber
	\end{align*}
\end{proposition}

\medskip

\begin{remark}\label{rem:obs_fns}
The subspace $\clC$ is identical to
the space of ``observable functions'' in~\cite[Lemma
9]{van2009observability}. An explanation for this correspondence is
provided 
in~\cite{kim2019observability} where it is shown that the BSDE is the
dual of the Zakai equation of filtering.  
\end{remark}

\medskip

%


Because $A$ is a stochastic matrix, its eigenvalues are either in the
open left half-plane or at zero.  To define stabilizability, consider
first the zero subspace:
\[
	S_0:=\{y \in \Re^d \mid \; Ay = 0\}
\]

\medskip


\begin{definition} 
	The BSDE~\eqref{eq:dual-bsde} is {\em stabilizable} if $S_0 \subset \clC$.
\end{definition}

\medskip

\begin{proposition}\label{prop:detectability-nullspace}
Consider the BSDE~\eqref{eq:dual-bsde}. Then
\begin{romannum}
\item If $\bS$ has a single ergodic class then BSDE is stabilizable.  
\item If $\bS=\cup_{k=1}^m \bS_k$ is partitioned into $m$ ergodic
  classes then the BSDE is stabilizable 
  if and only if the indicator functions $\ones_{\mathbb{S}_k} \in \clC$ for
  $k=1,2,\hdots,m$.  
\end{romannum}
\end{proposition}

\medskip


\section{Main result}\label{sec:main}

\begin{theorem}\label{thm:main-result}
Suppose the Wonham filter is stable. 
Then the BSDE~\eqref{eq:dual-bsde} is stabilizable.  Conversely,
if the BSDE is stabilizable then 
\[
\pi_T^{\nu}(\ones_{\bS_k}) \stackrel{(T\to\infty)}{\longrightarrow}
\ones_{\bS_k}(X_0) \quad \sP^{\mu}\text{-a.s.}
\] 
whenever $\mu \ll \nu$. (That is, the filter asymptotically
detects the correct ergodic class.)
\end{theorem}

\medskip

It is shown in Appendix~\ref{appdx:spliting-filters} that for any given $\nu \in \clP(\bS)$, one
can pick $\{\nu_k\in \clP(\bS):1\leq k \leq m\}$ such that $\nu_k$ has
support on $\bS_k$ and
\begin{equation}
\pi_T^\nu(f) = \sum_{k=1}^m \pi_T^\nu(\ones_{\bS_k})\pi_T^{\nu_k}(f) 
\label{eq:spliting-filters}
\end{equation}

Using~\Thm{thm:main-result}, the problem of filter
stability reduces to the problem of filter stability for each
ergodic class.  This is also the justification for treating the ergodic and
non-ergodic signal cases separately.  It is known that, for the type of
observations considered here, the filter ``inherits'' the stability
property of the Markov process~\cite[Theorem 4.2]{baxendale2004asymptotic}. The proof of this is far
from straightforward and spurred much research during the first
decade~\cite{van2010nonlinear,chigansky2009intrinsic,budhiraja2003asymptotic}.
Assuming it to be true, from~\Thm{thm:main-result} it follows that 
stabilizability is both a necessary and sufficient condition for filter stability.    

The dual optimal control approach of this paper is also useful to the study
of filter stability of ergodic signals (e.g., to obtain results on 
asymptotic convergence
of $\pi_T^{\nu_k}(f)$ above).  This is the subject of a companion
paper published in the proceedings of the conference~\cite{kim2021ergodic}.


\medskip

\begin{remark}
The sufficient condition stated in~\cite[Theorem
4.4]{baxendale2004asymptotic} correctly stress the importance of the
``identifying'' property of the filter to identify the correct ergodic
class~\cite[Lemma 6.3]{baxendale2004asymptotic}.  Subsequently, the
definition of detectability was first introduced
in~\cite{van2009observability,van2010nonlinear}.  For the Wonham
filter, the detectability property was shown to be equivalent to
filter stability~\cite[Theorem 2]{van2009observability}.    

Because of Remark~\ref{rem:obs_fns}, the stabilizability definition
taken together with the result in Theorem~1 represent dual
counterparts of these earlier definitions and results.  Additional
details on duality between controllability of the BSDE and
observability of the filter can be found in our prior paper
~\cite[Prop. 2]{kim2019observability}.  

It is also worthwhile to note that the appropriate observability and
detectability definitions were discovered  only after a decade of intense
research; c.f.,~\cite{van2010nonlinear,chigansky2009intrinsic}.  In
contrast, using duality these definitions are obtained quite naturally.

\end{remark}

\medskip

\begin{remark}
For the stability of the Kalman filter, the importance
of detectability is well known. Most proofs of filter stability rely on the analysis of the dual LQ optimal control
problem~\cite[Ch.~9]{xiong2008introduction},~\cite[Sec.~2]{ocone1996asymptotic}.
Stabilizability of the dual system is then {\em actually} the condition
that is used to obtain the proof of Kalman filter stability.  Our
results can thus be viewed as generalization of the Kalman
filter stability theory.  

\end{remark}


\section{Proof of the Main result}\label{sec:main_proof}

\subsection{Probability spaces}\label{ssec:prob_spaces}

Recall ${\sf P}^\mu$ is the probability measure indicative of the fact
that $X_0\sim \mu$.  For $\mu\ll \nu$, consider a probability
measure ${\sf P}^\nu$ on the
common measurable space $(\Omega,\clF)$ as in Sec.~\ref{ssec:model}. It is noted that 
\[
\frac{\ud \sP^\mu}{\ud \sP^{\nu}}(\omega) = \sum_{x\in\bS} \frac{\mu(x)}{\nu(x)}
\ones_{[X_0 = x]}(\omega)
\]  
Then $(X,Z)$ have the same transition law and if $X_0\sim \nu$ under
$\sP^\nu$ then  $X_0\sim \mu$ under $\sP^\mu$.  The expectation
under ${\sf P}^\nu$ is denoted ${\sf E}^\nu(\cdot)$.  The
solution of the Wonham filter $\pi_t^\mu =
\E^\mu(X_t\mid\clZ_t)$ and $\pi_t^\nu = \E^\nu(X_t\mid \clZ_t)$.

\medskip

In the settings of this paper, a filter is obtained by solving
the dual optimal control problem.  A user who (incorrectly) believes the prior
to be $\nu$ solves the optimal control problem under the
(incorrect) measure
$\sP^\nu$: 
\begin{align}\label{eq:opt_control_prob_Pnu}
\sJ_T^\nu (U;f)
= \E^\nu \Big( |Y_0(X_0)-\nu(Y_0)|^2 + \int_0^T \ell
  (Y_t,V_t,U_t\,;X_t) \ud t \Big)
\end{align}
subject to the BSDE constraint~\eqref{eq:dual-bsde}.  Note the two
changes: the expectation is now with respect to $\sP^\nu$ and $\nu(Y_0)$
appears in the terminal cost (first of the two terms). The optimal control for this problem
is denoted $U^\nu$ and the associated optimal trajectory is
$(Y^\nu,V^\nu)$.  The counterpart of~\eqref{eq:optimal-estimate-mu} is
\begin{equation}\label{eq:optimal-estimate-nu}
\pi_T^\nu(f) = \nu(Y_0^\nu) - \int_0^T U_t^{\nu} \ud Z_t\quad \sP^\nu-\text{a.s.}
\end{equation}
and $\sJ_T^\nu(f)  = \E^{\nu}( |f(X_T) - \pi_T^\nu(f) |^2 ) \leq
\frac{1}{4}|\text{osc}(f)|^2$ for all $T\geq 0$.  

\medskip

Before the estimator is assessed with respect to the (correct) measure
$\sP^\mu$, there are three technical concerns:
\begin{enumerate}
\item Existence and uniqueness of the optimal control $U^\nu$ and the solution $(Y^\nu,V^\nu)$ for all
  $T\geq 0$.  This ensures in particular that the
  righthand-side~\eqref{eq:optimal-estimate-nu} is well-defined.
\item Admissibility of the optimal control $U^\nu$ with respect to
  $\sP^\mu$.  This ensures that the optimal control can be assessed using $\sP^\mu$
  whenever $\mu\ll\nu$. 
\item Apriori bounds and continuity properties for the value
  $\sJ_T^\mu (U^\nu)$ as $\mu\to\nu$.   
\end{enumerate}

\medskip

Appendix~\ref{apdx:dual-optimal-control} contains the requisite
technical background that also serves to address these concerns.  


\subsection{Relationship of duality to filter stability}

For any $S_T\in L^2_{\clZ_T}$, the projection theorem 
gives
\begin{equation*}
\label{eq:pyth}
\E^{\mu} ( |f(X_T) - S_T |^2 ) =\sJ_T^\mu+ \E^{\mu} ( |\pi_T^\mu(f) - S_T |^2 )  
\end{equation*}
and using the duality formula~\eqref{eq:duality_reln} in the
left-hand side
\[
\E^{\mu} ( |\pi_T^\mu(f) - S_T |^2 ) = (\sJ_T^\mu (U) - \sJ_T^\mu) +
|\pi_0(Y_0) -\mu(Y_0)|^2
\]
With $\pi_0=\nu$ and $U=U^\nu$, the estimate
$S_T=\pi_T^\nu(f)$ and therefore
\[
\E^{\mu} ( |\pi_T^\mu(f) - \pi_T^\nu(f) |^2 ) = (\sJ_T^\mu (U^\nu) - \sJ_T^\mu) +
|\nu(Y_0^\nu) -\mu(Y_0^\nu)|^2
\]

Both the terms on the righthand-side are non-negative (the first term
so because
$\sJ_T^\mu$ is the minimum value and $U^\nu$ is $\sP^\mu$-admissible).  Therefore, the limit of the
lefthand-side (as $T\to\infty$) is $0$ if and only if each of the two terms on
the righthand-side individually approach $0$.  We state the result as
a proposition.    

\medskip

\begin{proposition}\label{prop:asymptotic-optimality}
The filter is stable if and only if 
	\begin{subequations}
		\begin{align}
		\mu(Y_0^\nu) - \nu(Y_0^\nu) \;\;
                  &\stackrel{(T\to\infty)}{\longrightarrow} \;\; 0 \label{eq:a1_st}\\
		\sJ_T^\mu(U^\nu)-\sJ_T^\mu
                  \;\;&\stackrel{(T\to\infty)}{\longrightarrow} \;\; 0 \label{eq:a2_st}
		\end{align}
	\end{subequations}
whenever $\mu\ll\nu$.
\end{proposition}

\medskip

\begin{remark}
These conditions are the nonlinear counterparts of the sufficient
conditions for the stability of the Kalman filter in~\cite[Theorem 2.3]{ocone1996asymptotic}.    
\begin{enumerate}
\item Equation~\eqref{eq:a1_st} means that the optimal control
  system is asymptotically stable.  That is, $Y_0\to \text{(const.)}
  \ones$ as $T\rightarrow\infty$.  This is also the reason why the stabilizability
  condition is important to the problem of filter stability.  The
  condition plays the same role in linear and nonlinear settings. 
\item Equation~\eqref{eq:a2_st} means that the value converges to the
  optimal value.  Since the optimal value $\sJ_T^\mu(f)$ has the
  interpretation of the minimum variance, its convergence is analogous to the
  convergence of the solution of the DRE in the Kalman
  filter.  In linear settings, the latter is deduced by establishing
  an asymptotic limit for the value function of the dual optimal control problem~\cite[Sec.~9.4]{xiong2008introduction}.   
\end{enumerate}



\end{remark}

\medskip

\subsection{Proof of necessity in \Thm{thm:main-result}}

In this subsection, we prove the necessity part of
\Thm{thm:main-result}:  That is, if the filter is stable then is the
BSDE~\eqref{eq:dual-bsde} is stabilizable.  The argument rests on the
result described in the following proposition whose proof appears in
the Appendix~\ref{apdx:pf-necessary-condition}.

\medskip

\begin{proposition}\label{prop:necessity}
Suppose the BSDE~\eqref{eq:dual-bsde} is not stabilizable.  Then
there exists an $f\in\clC^\perp$ such that for any $T$ and any
$U\in\clU$, the solution to the BSDE~\eqref{eq:dual-bsde} is given by 
\[
Y_0 = y_{0} + f 
\]
where $y_{0} \in \clC$ (and can depend upon $T$ and $U$). Since
$\ones\in\clC$, this implies $| Y_0 - \frac{1}{d} (\ones^\top Y_0) \ones
| \, \geq\, |f|$, i.e., $Y_0$ is uniformly bounded away from the
subspace of constant vectors.   

\end{proposition}

\medskip

The proof of necessity follows from \Prop{prop:necessity}:  Pick any $\nu\in
\clP(\bS)$ such that $0<\nu(x)<1$ for all $x\in\bS$ and set $\mu = \nu
+ \epsilon f$ where $\epsilon$ is chosen sufficiently small such that
$\mu\in\clP(\bS)$.  Then $|\mu(Y_0)-\nu(Y_0)| = \epsilon |f|^2$.  Applying
Prop.~\ref{prop:asymptotic-optimality}, the filter is not stable for
this choice of $\mu$ and $\nu$.  


\subsection{Completing the proof of  \Thm{thm:main-result}}

We first state a technical lemma that is used in the proof.  The proof
of the Lemma appears in Appendix~\ref{appdx:prop:stationary}. 

\medskip

\begin{lemma}\label{prop:stationary}
Suppose $\barpi$ is an invariant measure of $A$ (i.e., $A^\top \bar\mu  =
0$).  Then for each fixed $f\in\Re^d$
\begin{romannum}
\item The sequence $\{\sJ_T^{\bar{\mu}}(f):T\geq 0\}$ is bounded,
  non-negative, and non-increasing in $T$.  Therefore,
  $\sJ_T^{\bar{\mu}}(f)$ converges as $T\rightarrow \infty$.  Denote the
  limit as $\sJ_\infty^{\bar{\mu}}(f)$.
\item For a given $\mu\in{\cal P}(\mathbb{S})$, denote $\mu_T :=
  e^{A^\top T} \mu$.  Suppose $\mu_T\to \barpi$ as $T\to\infty$.  Then
\[
\limsup_{T\to\infty} \sJ_{T}^\mu(f) \le \sJ_{\infty}^{\barpi}(f)
\]
\end{romannum}
\end{lemma}

\medskip

We now complete the proof of \Thm{thm:main-result}: That is, we show
that if $\ones_{\bS_k}\in\clC$ then
\[
\pi_T^{\nu}(\ones_{\bS_k}) \stackrel{(T\to\infty)}{\longrightarrow}
\ones_{\bS_k} (X_0)\quad \sP^{\mu}\text{-a.s.}
\] 
The proof is in the following three steps:
\medskip

\noindent 1. In step~1, we show that 
$\pi_T^\nu(\ones_{\bS_k})$ converges ${\sf P}^\nu$-a.s.  

\medskip

\noindent 1. In step~2, we show that if 
$\ones_{\bS_k}\in\clC$ then $\sJ_T^\barpi(\ones_{\bS_k}) \to 0$ as $T\to\infty$ where
$\barpi$ is any invariant measure of $A$.  We use part (i) of
\Lemma{prop:stationary} to prove this result.  
\medskip

\noindent 3. In step~3, we combine the conclusions of steps~1 and~2 to prove
  the result.  We use part (ii) of
\Lemma{prop:stationary} to prove this result.

\newP{Step 1} Consider the Wonham filter~\eqref{eq:Wonham} with $\pi_0=\nu$.  Since
$A \ones_{\bS_k} = 0$, $\{\pi_T^\nu(\ones_{\bS_k}):T\geq 0\}$ is a
bounded ${\sf
  P}^\nu$-martingale and therefore converges $\sP^\nu$-a.s.
(Therefore, the a.s. convergence does not require stabilizability of the
model.)

\newP{Step 2} Suppose $\barpi$ is any invariant measure.  Then $\sJ_T^\barpi$
is monotone (part~(i) of Lemma~\ref{prop:stationary}).  In the
following, we construct a
sequence of admissible control input $\{U^{(T)}:T=1,2,\hdots\}$ such that 
$\sJ_T^\barpi(U^{(T)};\ones_{\bS_k}) \to 0$ as $T\to\infty$.  Since
$\sJ_T^\barpi(\ones_{\bS_k})$ is the minimum value this implies
$\sJ_T^\barpi(\ones_{\bS_k}) \to 0$ as $T\to\infty$ (for this
particular sub-sequence).  Since $\sJ_T^\barpi$
is monotone, the limit exists and equals this sub-sequential
limit.

Suppose $\ones_{\bS_k} \in \clC$.  Then we claim that there exists an admissible
control $U^{(1)}:=\{U_t^{(1)}:0\leq t \leq 1\}$ and a constant
$c\in\Re$ such that $Y_T^{(1)} =
\ones_{\bS_k}$ and $Y_0^{(1)} = c \ones$.  The claim follows from the
definition of controllable space.  A proof of the claim appears at the
end of the proof.   
Assuming the claim to be true for now, denote the associated
solution of the BSDE~\eqref{eq:dual-bsde} as
$(Y^{(1)},V^{(1)}):=\{(Y_t^{(1)},V_t^{(1)}): 0\leq t \leq 1\}$.  

Since $Z$ is a w.p. under the Girsanov change of measure, there exists
a functional $\phi(\cdot\,,\,\cdot):[0,1]\times C([0,1];\Re)\to\Re$ such
that 
\[
U_t^{(1)} = \phi(t\,,\,\{Z_s:0\leq s \leq t\}) \quad \sP^{\bar\mu}-\text{a.s.}
\]

Now consider the following control over the time-horizon $[0,n]$:  For $l=0,1,2,\hdots,n-1$
\[
U_t^{(n)} := \frac{1}{n}\phi\big(t-l\,,\, \{Z_s:l\leq s \leq t\}\big),\quad t\in(l,(l+1)]
\]  
Such a control input is clearly admissible.  
With $Y_n = \ones_{\bS_k}$, one obtains the following solution
$(Y^{(n)},V^{(n)}):=\{(Y_t^{(n)},V_t^{(n)}): 0\leq t \leq n\}$ of the
BSDE~\eqref{eq:dual-bsde}:  
\[
V_t^{(n)} \stackrel{\text{(d)}}{=} \frac{1}{n}V_{t-l}^{(1)},\quad t\in(l,(l+1)]
\]
for $l=0,1,2,\hdots,n-1$, $V_0^{(n)} = V_0^{(1)}$, and
\begin{align*}
Y_t^{(n)} \stackrel{\text{(d)}}{=} \begin{cases} \frac{1}{n}Y_t^{(1)}
  + \frac{n-1}{n} \ones_{\bS_k} &\mbox{if } t \in ((n-1),n] \\
\frac{1}{n}Y_t^{(1)} + \frac{n-2}{n}\ones_{\bS_k} + \frac{c}{n}\ones &
\mbox{if }t \in ((n-2),(n-1)] \\
\quad\vdots & \quad\vdots \\
\frac{1}{n}Y_t^{(1)} +\frac{1}{n}\ones_{\bS_k} + \frac{(n-2)c}{n}\ones & \mbox{if }t \in (1,2] \\
\frac{1}{n}Y_t^{(1)} +
\frac{(n-1)c}{n}\ones & \mbox{if }t \in (0,1] 
\end{cases}
\end{align*}
and $Y_0^{(n)} = c \ones$.  

Since $Y_0^{(n)} = c\ones$, the terminal cost 
$|Y_0^{(n)} (X_0)-\bar\mu(Y_0^{(n)})|^2 = 0$. And since $X_t\sim \bar\mu$, for $l=0,1,2,\hdots,n-1$:
\begin{align*}
\E^{\barpi}\Big(&\int_{l}^{(l+1)}  \Gamma(Y_t^{(n)})(X_t)  + 
                  |U_t^{(n)} + V_t^{(n)} (X_t)|_R^2 \ud t \Big) \\
&=\frac{1}{n^2}\E^{\barpi}\Big(\int_{0}^{1} \Gamma(Y_t^{(1)})(X_t) +
  |U_t^{(1)}+ V_t^{(1)} (X_t)|_R^2 \ud t \Big)
\end{align*}
Therefore,
\[
\sJ_{n}^{\barpi}(U^{(n)}) =\frac{1}{n}\sJ_1^{\barpi}(U^{(1)})
\]
and thus the optimal value 
\begin{equation}\label{eq:l2-conv-indicator}
\sJ_n^\barpi \le\sJ_{n}^{\barpi}(U^{(n)}) =\frac{1}{n}\sJ_1^{\barpi}(U^{(1)}) \longrightarrow 0 \quad \text{as}\quad n\to\infty
\end{equation}

\newP{Step 3} Suppose $\nu \in\clP(\bS)$ and $\ones_{\bS_k} \in
\clC$.  In this final step, we show that $\sJ_T^\nu \to 0$ and 
\[
\pi_T^{\nu}(\ones_{\bS_k}) \stackrel{(T\to\infty)}{\longrightarrow}
\ones_{\bS_k}(X_0) \quad \sP^{\nu}\text{-a.s.}
\] 

Let $\barpi_l\in\clP(\bS)$ be the invariant measure for the
$l^{\text{th}}$-ergodic class and $a_l := \nu (\ones_{\bS_l})$ for
$l=1,2,\hdots,m$.  Choose the invariant measure as follows:
\[
\barpi = a_1\barpi_1 + a_2\barpi_2 + \hdots + a_m\barpi_m 
\]
From step~2, we know that $\sJ_\infty^\barpi = 0$.  Also, $\nu_T :=
e^{A^\top T}\nu\to \barpi$ as $T\to\infty$. 
Therefore, using part~(ii) of Lemma~\ref{prop:stationary}, 
\[
\limsup_{T\to \infty} \sJ_T^\nu \leq \sJ_\infty^\barpi  = 0
\]
which shows that $\sJ_T^\nu \to 0$ as $T\to\infty$.   

Since $\bS_k$ is an ergodic class,
\[
\sJ_T^\nu = \E^\nu\big(|\ones_{\bS_k}(X_T)  -
\pi_T^\nu(\ones_{\bS_k})|^2\big) = \E^\nu\big(|\ones_{\bS_k}(X_0)  - \pi_T^\nu(\ones_{\bS_k})|^2\big)
\]
By Fatou's lemma,
\[
\E^\nu\big(\liminf_{T\to\infty}|\pi_T^\nu(\ones_{\bS_k})-\ones_{\bS_k}(X_0)|^2 \big) \le \lim_{T\to \infty}\sJ_T^\nu = 0
\]
In step 1, we showed that $\pi_T^\nu(\ones_{\bS_k})$ has an
a.s. limit.  So, $\liminf$ is replaced as
\[
\lim_{T\to \infty} |\pi_T^\nu(\ones_{\bS_k})-\ones_{[X_0\in\bS_k]}|^2 = 0\quad \sP^\nu\text{-a.s.}
\]
and therefore also $\sP^\mu$-a.s. whenever $\mu\ll \nu$.

\medskip

\newP{Proof of the claim in step~2} Suppose $f \in \clC_T$. Since
$\clC_T$ is $A$-invariant (see Prop.~\ref{thm:rank}), $e^{AT}f \in
\clC$. Therefore, from definition of $\clC$, there is a deterministic constant
$c\in\Re$ and an admissible control $U\in\clU$ such that the solution of the
BSDE~\eqref{eq:dual-bsde} is obtained with $Y_T = c\ones$ and $Y_0 =
e^{AT}f$. Now consider a second solution of the BSDE~\eqref{eq:dual-bsde} with $Y_T =
c\ones+f$ and zero control input. Since $A\ones = 0$, this second solution
is $(Y_t,V_t) = (c\ones + e^{A(T-t)}f, 0)$ for $t \in [0,T]$. By
linearity, we subtract the two solutions to show that with $Y_T = f$
and control $-U$, one obtains $Y_0 = c\ones$. \qed

 \section{Conclusions and Directions for future work}

Nonlinear filtering is an old subject. It is also notoriously
difficult, which is why it remains an exciting research domain with
many open questions remaining.    This paper presents a new attack on
filter stability, which we hope will open new avenues for research on
filter performance and design. There are several avenues for future work:

Although we prove a filter stability result when the dual model is
stabilizable, a more nuanced understanding is possible through the
consideration of the controllable subspace ${\cal C}$ of the
BSDE. Specifically, ${\cal C}$ is the space of functions for
which the filter forgets the initial condition.  A duality-based proof
of this remains open.

Another important question is to relate our work to the deterministic
  definitions. Given the importance of the IOSS definitions of
observability for deterministic models, it is of interest to
investigate the dual optimal control
problem in its small noise limit (in
particular as $R\downarrow 0$).  




\bibliographystyle{IEEEtran}
\bibliography{duality,backward_sde,filter-stability-observability,jin_papers}

\appendix

\section{Appendix}

\subsection{Technical background} \label{apdx:dual-optimal-control}

In this section, we describe some technical background on well-posedness
of the optimal control problem~\eqref{eq:opt_control_prob_Pnu} with
respect to $\sP^\nu$.  We assume $\mu,\nu\in\sP(\bS)$ and $\mu\ll\nu$.  

Define the innovation process $I^\nu:= \{I_t^\nu\in\Re:t\in [0,T]\}$ and the
covariance process $\Sigma^\nu:=\{\Sigma_t^\nu \in\Re^{d\times d} :
0\le t\le T\}$ as follows:
\begin{align*}
\text{(innovation)}:\;\quad I_t^\nu &:= Z_t - \int_0^t \pi_s^\nu(h) \ud s\\
\text{(covariance)}:\quad \Sigma_t^\nu &:= \dv
(\pi_t^\nu)-\pi_t^\nu(\pi_t^{\nu})^\top
\end{align*}
It is known that $I^\nu$ is $\sP^\nu$-w.p. and moreover the
filtration generated by $Z$ and $I$ are identical~\cite{allinger1981new}.


The following identity is proved in~\cite[Theorem
5]{kim2019duality}: For any admissible control input $U$
\begin{align*}
\sJ_T^\nu(U) &= \E^\nu \Big(Y_t^\top \Sigma_t^\nu Y_t + \int_t^T \ell (Y_s,V_s,U_s;X_s)\ud s\Big)\\
&\quad \quad +\E^\nu \Big(\int_0^t \big|U_s+ (h^\top \Sigma_s^\nu Y_s+\pi_s^\nu(V_s))\big|^2\ud s\Big)
\end{align*}
From this identity, the following conclusions are deduced:
\begin{romannum}
\item The minimum value $\sJ_T^\nu = \E^\nu (f^\top \Sigma_T^\nu f) $.  
\item The optimal control is of the feedback form
\begin{equation}\label{eq:opt_control_law}
U_t^\nu = -h^\top \Sigma_t^\nu Y_t - \pi_t^\nu(V_t), \quad \leq t\leq T
\end{equation}
\end{romannum}

\medskip

Using the optimal control law~\eqref{eq:opt_control_law} in the
BSDE~\eqref{eq:dual-bsde} results in the following linear feedback control
system:
\begin{align}
-\ud Y_t &= \big((A - h h^\top \Sigma_t^\nu)Y_t + h\cdot V_t -
h\pi_t^\nu(V_t) - V_t \pi_t^\nu(h) \big)\ud t - V_t \ud I_t^\nu \nonumber \\
Y_T &= f 
\label{eq:dual-bsde_opt}
\end{align} 

The following proposition provides the answer to the three concerns in~\Sec{ssec:prob_spaces}:

\medskip

\begin{proposition}
Consider the linear BSDE~\eqref{eq:dual-bsde_opt} with $f\in\Re^d$ and
$T> 0$.
Then
\begin{romannum}
\item There exists a unique solution $(Y^\nu,V^\nu) \in
  L^2_{\clZ}([0,T];\Re^d\times \Re^d)$.  The optimal
  value $\sJ_T^\nu(U^\nu;f) \le \frac{1}{4}
  |\text{osc}(f)|^2$. 
\item Suppose $\mu\ll \nu$.  The optimal control $U^\nu$ is in $L^2_{\clZ}([0,T];\Re)$ also with
  respect to the $\sP^\mu$-measure (so it is admissible).  The value
\begin{align*}
\sJ_T^\mu(U^\nu;f) + |\mu(Y_0^\nu) - \nu(Y_0^\nu)|^2 \le \max_{x\in\bS}\frac{\mu(x)}{\nu(x)} \frac{1}{4}
  |\text{osc}(f)|^2
\end{align*}
\item (Continuity property) Consider a family $\{\nu_n\in\sP(\bS) :  n
  =1,2,\hdots\}$ such that $\nu_n\ll \nu$ and $\nu_n \to \nu$.  Then
\begin{equation}\label{eq:continuity-property}
|\sJ_T^{\nu_n}(U^\nu;f) - \sJ_T^\nu(U^\nu;f)| \longrightarrow 0 \tag{CP}
\end{equation}
where the convergence is uniform in $T$.  

\end{romannum}
\end{proposition}

\medskip

{\it Proof:}

\newP{Part (i)} The optimal control system~\eqref{eq:dual-bsde_opt} is
a linear BSDE with random but bounded coefficients.  The coefficients
are bounded because each element of $\pi_t^\nu$ and $\Sigma_t^\nu$ is in $[0,1]$.  Part~(i) follows
from existence uniqueness theory of linear BSDE~\cite[Theorem 7.2.2]{yong1999stochastic}.  Because
$\sJ_T^\nu = \E^\nu (f^\top \Sigma_T^\nu f) $, the uniform bound
readily follows.  

\newP{Part (ii)}  Since $\mu \ll \nu$, $\sP^\mu \ll \sP^\nu$ with
$\frac{\ud \sP^\mu}{\ud \sP^\nu} =
\frac{\mu(X_0)}{\nu(X_0)}$. So,
\begin{align*}
\E^\mu\big(\int_0^T U_t^2 \ud t\big) &=
  \E^\nu\big(\frac{\mu(X_0)}{\nu(X_0)} \int_0^T U_t^2 \ud t\big) \\ 
& \le \max_{x\in\bS}\frac{\mu(x)}{\nu(x)} \E^\nu\big(\int_0^T U_t^2 \ud t\big)
\end{align*}
Therefore if $U$ is $\sP^\nu$-admissible then it is also
$\sP^\mu$-admissible.  From duality~\eqref{eq:duality_reln},
\begin{align*}
\sJ_T^\mu(U^\nu;f) + |\mu&(Y_0^\nu) - \nu(Y_0^\nu)|^2 = \E^\mu\big(|f(X_T)-\pi_T^\nu(f)|^2\big)
\end{align*}
and using the change of measure the righthand-side
\[
\E^\mu\big(|f(X_T)-\pi_T^\nu(f)|^2\big) \leq \max_{x\in\bS}\frac{\mu(x)}{\nu(x)}\E^\nu\big(|f(X_T)-\pi_T^\nu(f)|^2\big)
\]
The result follows from using the bound from part (i).

\newP{Part (iii)} We have
\begin{align*} \sJ_{T}^{\nu_n}(U^\nu;f) 
= 
 \E^{\nu_n} \Big( |Y_0^\nu(X_0) & -\nu_n(Y_0^\nu)|^2 \\ & + \int_0^{T}
\ell(Y_t^\nu,V_t^\nu,U_t^\nu;X_t)\ud t \Big)
\end{align*}
We show $\sJ_{T}^{\nu_n}(U) \longrightarrow \sJ_{T}^\nu (U)$ if
$\nu_n\to \nu$. We consider each of the two terms:
\begin{enumerate}
\item The first term is written as
  $Y_0^{\nu\top}\Sigma_0^{\nu_{n}}Y_0^\nu$.  We have
\begin{align*}
\|\Sigma_0^{\nu_{n}} - \Sigma_0^{\nu}\|_2 &= \|\dv(\nu_{n}-\nu) - (\nu_{n}\nu_{n}^\top - \nu\nu)\|_2\\
&\le \|\dv(\nu_{n}-\nu)\|_2 + \|\nu_{n}\nu_{n}^\top - \nu\nu^\top\|_2\\
&\le 3\|\nu_{n}-\nu\|_{\infty} \;\;\stackrel{(n\to\infty)}{\longrightarrow} \;\; 0
\end{align*}  
By part (i) $Y_0^\nu
\Sigma_0^{\nu}Y_0^\nu$ is uniformly bounded so the limit is well-defined.  

\item
For the integral term, let $\xi_T :=\int_0^{T}
\ell(Y_t^\nu,V_t^\nu,U_t^\nu;X_t)\ud t$. Then
\begin{align*}
\E^{\nu}(\xi_T) &= \sum_{i=1}^d \nu_{n}(i)\E^{\delta_i}(\xi_T),\;\; 
\E^{\nu_{n}}(\xi_T) = \sum_{i=1}^d \nu_{n}(i)\E^{\delta_i}(\xi_T)
\end{align*} 
Since $\E^{\nu}(\xi_T) \leq \sJ_T^\nu$,
\[
|\E^{\nu_{n}}(\xi_T) - \E^\nu(\xi_T)| = \sum_{i=1}^d |\nu_{n}(i)-\nu(i)|\E^{\delta_i}(\xi_T) \;\; \longrightarrow \;\; 0
\]
and the convergence is uniform in $T$.  
\qed
\end{enumerate}

\medskip

\subsection{Proof of Proposition~\ref{prop:detectability-nullspace}} \label{apdx:prop:detectability-nullspace}
 
%
%
\begin{romannum}
\item For the ergodic case, $0$ is a simple eigenvalue of the matrix
  $A$ and $S_0=\text{span}\{\ones\}$.  
\item  Suppose $\bS=\cup_{k=1}^m \bS_k$ is a ergodic partition.  By
  choosing an appropriate coordinate, the rate matrix
\[
 A = \begin{pmatrix}
 A_1 & 0 & \cdots & 0\\
 0 & A_2 & \cdots & 0\\
 \vdots & \vdots & \ddots & \vdots\\
 0 & 0 & \cdots & A_m
 \end{pmatrix}
\]
where $A_k$ has a simple eigenvalue at $0$ with the associated eigenvector
$\ones_{\mathbb{S}_k}$.  (This is so because $\bS_k$ is an ergodic class.) Therefore,
\[
S_0 = \text{span}\{\ones_{\mathbb{S}_1}, \ones_{\mathbb{S}_2},\hdots, \ones_{\mathbb{S}_m}\}
\]
\end{romannum}

\subsection{Proof of the splitting~\eqref{eq:spliting-filters}}\label{appdx:spliting-filters}

Suppose $\mathbb{S} = \cup_{k=1}^m \mathbb{S}_k$ is an ergodic
partition.   
For each such ergodic class with
$\sP^\nu([X_0\in\bS_k]) > 0$ define 
\[
\nu_k(x) := \left\{ \begin{array}{cc}
                      \dfrac{\nu(i)}{\sP^\nu([X_0\in\bS_k])} &
                      \text{if}\;\; x\in \bS_k \\ 0 &
                      \text{if}\;\; x\notin \bS_k \end{array} \right.
\]
Clearly $\nu_k\ll\nu$ and 
\[
\frac{\ud \sP^{\nu_k}}{\ud \sP^\nu} (\omega) = \sum_{x\in\bS} \frac{\nu_k(x)}{\nu(x)}\ones_{[X_0=x]}(\omega) = \frac{\ones_{[X_0(\omega)\in\bS_k]}(\omega)}{\sP^\nu([X_0\in\bS_k])}
\]
An application of the Bayes' formula gives
\[
\E^{\nu_k}\big(f(X_T)|\clZ_T\big) = \frac{\E^\nu\big(f(X_T)\ones_{[X_0\in\bS_k]}|\clZ_T\big)}{\E^\nu\big(\ones_{[X_0\in\bS_k]}|\clZ_T\big)}
\]
and therefore
\begin{equation}\label{eq:Bayes_Identity}
\E^\nu\big(f(X_T)\ones_{[X_0\in\bS_k]}|\clZ_T\big) = \pi_T^\nu(\ones_{\bS_k})\pi_T^{\nu_k}(f) 
\end{equation}
where we have used the fact that
\[
\ones_{[X_0(\omega) \in\bS_k]} (\omega) =
\ones_{[X_T(\omega) \in\bS_k]}(\omega) \quad {\sf P}^\nu-a.s.
\]
Note that the identity~\eqref{eq:Bayes_Identity} is true for all
$k=1,2,\hdots,m$.  (If $\sP^\nu([X_0\in\bS_k]) = 0$ then both sides
are zero.)  Upon summing the identity over the index $k$, one arrives at
\[
\pi_T^\nu(f) = \sum_{k=1}^m \pi_T^\nu(\ones_{\bS_k})\pi_T^{\nu_k}(f) 
\]
If $\sP^\nu([X_0\in\bS_k]) = 0$ then take $\nu_k$ to be any 
probability measure with support on $\bS_k$ (e.g., the invariant measure for the restriction of the Markov process on $\bS_k$).  

\subsection{Proof of \Prop{prop:necessity}}\label{apdx:pf-necessary-condition}

Suppose $\clC$ is the controllable subspace with dimension strictly
less than $d$. Consider  the splitting $\Re^d = \clC \oplus
\clC^\perp$ and an associated orthogonal transformation $\tf: \clC \oplus \clC^\perp
\rightarrow \Re^d$ such that 
\[
Y_t = \tf\bar{Y}_t,\; V_t = \tf\bar{V}_t \quad \text{where} \quad \bar{Y}_t =
\begin{bmatrix} \bar{Y}_t^{c} \\ \bar{Y}_t^{uc} \end{bmatrix},\; 
 \bar{V}_t =
\begin{bmatrix} \bar{V}_t^{c} \\ \bar{V}_t^{uc} \end{bmatrix}
\] 
With respect to the new coordinates, the BSDE~\eqref{eq:dual-bsde} becomes
\begin{align*}
-\ud \bar{Y}_t = \big(\bar{A}\bar{Y}_t + \bar{h} U_t +
  \bar{K}\bar{V}_t\big) \ud t - \bar{V}_t \ud Z_t,\quad \bar{Y}_T =
  \tf^{-1}f =\bar{f} = \begin{bmatrix} \bar{f}^{c} \\ \bar{f}^{uc} \end{bmatrix}
\end{align*}
where (because $\clC$ is $A$- and $\dv(h)$-invariant) the matrices
\[
\bar{A} = \tf^{-1}A\tf = \begin{bmatrix} \bar{A}_{c} & * \\ 0 &
  \bar{A}_{uc} \end{bmatrix}, \quad \bar{K} = \tf^{-1}\dv(h) \tf = \begin{bmatrix}
\bar{K}_{c} & * \\ 0 & \bar{K}_{uc}
\end{bmatrix}
\] 
and $\bar{h}=\tf^{-1} h=\begin{bmatrix}
\bar{h}_c  \\ 0
\end{bmatrix}$. 
On the $\clC^{\perp}$ subspace
\[
-\ud \bar{Y}_t^{uc} = \big(\bar{A}_{uc} \bar{Y}_t^{uc} + K_{uc}
\bar{V}_t^{uc}\big) \ud t - \bar{V}_t^{uc} \ud Z_t,\quad
\bar{Y}_T^{uc} = \bar{f}^{uc}
\]
whose unique solution is given by
\begin{equation}\label{eq:C-complement}
\bar{Y}_t^{uc} = e^{\bar{A}_{uc} (T-t)}\bar{f}^{uc},\quad
\bar{V}_t^{uc} \equiv 0,\quad t\in[0,T]
\end{equation}
This is so because $f$ is a deterministic function.  
  
If the BSDE is not stabilizable, then there exists a non-zero
vector $\bar{\eta}$ such that $\bar{A}_{uc}\bar\eta = 0$.  Set  
\[
\bar{f} =\begin{bmatrix} 0 \\ \bar{\eta} \end{bmatrix} \quad
\Longrightarrow \quad \bar{Y}_0 = \begin{bmatrix} * \\
  \bar{\eta} \end{bmatrix}
\]   
Since $\ones \in \clC$, the length of the vector $\left. {Y}_0
\right|_{\ones^\perp}$ is at least $|\bar{\eta}|=|\bar{f}|=|f|$.

\subsection{Proof of Lemma~\ref{prop:stationary}}\label{appdx:prop:stationary}

The proof of the lemma requires a technical construction. Consider the time horizon $[0,T_1+T_2]$.  If $X_0\sim \mu$ then
$X_{T_1}\sim e^{A^\top T_1}\mu=:\mu_{T_1}$.  This is useful to relate
the properties of $\sfJ^\mu_{T_1+T_2}(\cdot)$ and
$\sfJ^{\mu_{T_1}}_{T_2}(\cdot)$.  For this purpose, consider first the
time horizon $[T_1,T_1+T_2]$.  Over this time horizon, introduce the filtration 
\[
\tilde\clZ_{t-T_1}:=\{Z_t - Z_{T_1}\;:\; T_1 \le  t \le T_1+T_2\}
\]
For a control $\tilde{U} \in L^2_{\tilde{\clZ}}([0,T_2])$, let
$\{(\tilde{Y}_t, \tilde{V}_t): t\in[0,T_2]\}$ denote the solution of the
BSDE~\eqref{eq:dual-bsde} with $\tilde{Y}_{T_2}=f$.  The control $\tilde{U}$ is extended to the
time-horizon $[0,T_1+T_2]$ as follows:
\begin{equation}\label{eq:zero-U-control}
U_t=\begin{cases}
0 & 0 \le t < T_1\\
\tilde{U}_{t-T_1} & T_1 \le t \le T_1+T_2
\end{cases}
\end{equation}
The control $U\in{\cal U}$ and yields the following solution of
the BSDE~\eqref{eq:dual-bsde}:
\[ \label{eq:soln-barU}
(Y_t,V_t)=\begin{cases}
(e^{A(T_1-t)}\tilde{Y}_{0},\;0) & 0 \le t < T_1\\
(\tilde{Y}_{t-T_1},\;\tilde{V}_{t-T_1}) & T_1 \le t \le T_1+T_2
\end{cases}
\]
Under this definition, we claim that  
\begin{equation}\label{eq:claim1}
\sfJ^\mu_{T_1+T_2}(U) = \sfJ^{\mu_{T_1}}_{T_2}(\tilde{U})
\end{equation}
and then the two results in the lemma are direct consequences of the claim:

\newP{Part (i)} Take $\mu=\barpi$
and  $\tilde{U} = U^{\barpi}$.  Then
\[
\sJ_{T_1+T_2}^{\barpi} \le \sJ_{T_1+T_2}^{\barpi}(U)
\overset{\text{Eq.~\eqref{eq:claim1}}}{=\joinrel=} \sJ_{T_2}^{\barpi}(U^{\barpi}) = \sJ_{T_2}^{\barpi}
\] 
where we used the facts that (i) $\mu_{T_1}=e^{A^\top T_1}\barpi =
\barpi$ because $\barpi$ is the invariant measure; and (ii)
$U^{\barpi}$ is the optimal control for the $\sJ_{T_2}^{\barpi}(\cdot)$
problem.  Therefore, $\sJ_{T}^{\barpi}$ is monotone in $T$ 
and converges as $T\rightarrow \infty$.  Denote the
  limit as $\sJ_\infty^{\bar{\mu}}$.

\medskip

\newP{Part (ii)} Let $T=T_1+T_2$.  For $\mu\in{\cal P}(\mathbb{S})$,
with $\tilde{U} = U^{\barpi}$
\begin{equation}\label{eq:cor1_part2}
\sJ_T^\mu \le \sJ_{T_1+T_2}^{\mu}(U)
\overset{\text{Eq.~\eqref{eq:claim1}}}{=\joinrel=} \sJ_{T_2}^{\mu_{T_1}}(U^\barpi) 
\end{equation}
We have
\[
|\sJ_{T_2}^{\mu_{T_1}}(U^{\barpi}) - \sJ_\infty^{\bar{\mu}}| \leq
|\sJ_{T_2}^{\mu_{T_1}}(U^{\barpi}) -
\sJ_{T_2}^{\barpi}(U^{\barpi})| + |\sJ_{T_2}^{\barpi}(U^{\barpi}) - \sJ_\infty^{\bar{\mu}}|
\]
The second term on the righthand-side does not depend upon
$T_1$. Because $U^{\barpi}$ is the optimal control input, this term
goes to zero as $T_2\to \infty$:  That is, given $\epsilon>0$, there
exist an $n_2$ such that 
\[
|\sJ_{T_2}^{\barpi}(U^{\barpi}) -
\sJ_\infty^{\bar{\mu}}|\leq \epsilon\quad \forall T_2 \ge n_2
\]
Now fix $T_2 = n_2$ and apply continuity property~\eqref{eq:continuity-property} to the first
term on the righthand-side: There exists $n_1=n_1(n_2)$ such that 
\[
|\sJ_{n_2}^{\mu_{T_1}}(U^{\barpi}) -
\sJ_{n_2}^{\barpi}(U^{\barpi})| \leq \epsilon \quad \forall T_1 \ge n_1
\]
Combine these inequalities concludes for all $T_1 \ge n_1$,
\[
\sJ_{n_2}^{\mu_{T_1}}(U^\barpi) \leq
\sJ_\infty^{\bar{\mu}} + 2\epsilon
\]    
From~\eqref{eq:cor1_part2}, $
\sJ_T^\mu \le \sJ_{n_2}^{\mu_{T-n_2}}(U^\barpi)$. Therefore, for all $T \ge n_1+n_2$, 
\[
\sJ_{T}^\mu \le \sJ_{n_2}^{\mu_{T-n_2}}(U^\barpi) \leq \sJ_\infty^{\bar{\mu}} + 2\epsilon
\]
Since $\epsilon$ is arbitrary, the result follows.  

\medskip

It remains to prove the claim~\eqref{eq:claim1}.
We have
\begin{align*}
\sJ_{T}^\mu(U) &= \E^{\mu}\Big(|Y_0(X_0)-\mu(Y_0)|^2 +\int_0^{T_1} \Gamma (Y_t)(X_t) \ud
t\Big) \\
&\quad+\quad \E^{\mu} \Big(\int_{T_1}^{T_1+T_2} \Gamma (Y_t)(X_t)
+ |U_t+V_t(X_t)|_R^2\ud t \Big)
\end{align*}
Each of the two terms is simplified separately. 

Consider the control input $U$ defined according
to~\eqref{eq:zero-U-control}.  
Since $Y_{T_1} = \tilde{Y}_{0}$ is a deterministic function and the control is set to be zero, $V_t = 0$ on $0\le t < T_1$ and the BSDE becomes ODE:
\[
-\frac{\ud}{\ud t} Y_t = AY_t,\quad Y_{T_1} = \tilde{Y}_0
\]
A straightforward calculation shows that
\begin{align*}
&\E^{\mu}\Big(|Y_0(X_0)-\mu(Y_0)|^2 +\int_0^{T_1} \Gamma (Y_t)(X_t) \ud
t\Big) \\
&= \E^\mu\Big(|Y_{T_1}(X_{T_1})-\mu_{T_1}(Y_{T_1})|^2\Big)= \E^{\mu_{T_1}}\Big(|\tilde{Y}_{0}(X_0)-\mu_{T_1}(\tilde{Y}_{0})|^2\Big)
\end{align*}

The second term
\begin{align*}
& \E^{\mu} \Big(\int_{T_1}^{T_1+T_2}  \Gamma (Y_t)(X_t)  +
  |U_t+V_t(X_t)|_R^2\ud t \Big) \\
& = \E^{\mu} \Big(\int_{T_1}^{T_1+T_2}  \Gamma (\tilde{Y}_{t-T_1})(X_t) +  |\tilde{U}_{t-T_1}+\tilde{V}_{t-T_1}(
  X_t)|_R^2\ud t \Big)\\
& = \E^{\mu_{T_1}} \Big(\int_{0}^{T_2}  \Gamma(\tilde{Y}_{t})({X}_{t})+  |\tilde{U}_{t}+\tilde{V}_{t}( {X}_{t})|_R^2\ud t \Big)
\end{align*}
Combining the results of the two calculations yields:
\begin{align*}
\sJ_{T}^\mu(U)
&= \E^{\mu_{T_1}}\Big(|\tilde{Y}_{0}(X_0)-\mu_{T_1}(\tilde{Y}_{0})|^2\Big)\\
&\quad\quad+ \E^{\mu_{T_1}}\Big(\int_{0}^{T_2}  \Gamma(\tilde{Y}_{t})(\tilde{X}_{t}) +  |\tilde{U}_{t}+\tilde{V}_{t}( \tilde{X}_{t})|_R^2\ud t \Big)\\
&= \sJ_{T_2}^{\mu_{T_1}}(\tilde{U}) 
\end{align*}

\end{document}